\documentclass[12pt]{article}

\usepackage{amsmath, epsfig, cite}
\usepackage{amssymb}
\usepackage{amsfonts}
\usepackage{latexsym}

\newtheorem{thm}{Theorem}[section]

\newtheorem{conj}[thm]{Conjecture}
\newtheorem{lem}[thm]{Lemma}
\newtheorem{prob}[thm]{\it Problem}
\newcommand{\pf}{\noindent{\it Proof.} }
\newcommand*{\fl}[2]{\left\lfloor\frac{#1}{#2}\right\rfloor}

\def\fl#1{\left\lfloor#1\right\rfloor}

\numberwithin{equation}{section}

\newcommand{\qed}{{\hfill$\square$}\medskip}

\begin{document}

%\linenumbers
\nocite{*}
\begin{center}
{\Large\bf Proof of Sun's conjectures on integer-valued\\[5pt] polynomials}
\end{center}

\vskip 2mm \centerline{Victor J. W. Guo}
\begin{center}
{School of Mathematical Sciences, Huaiyin Normal University, Huai'an, Jiangsu 223300,
 People's Republic of China\\
{\tt jwguo@hytc.edu.cn} }

\end{center}

%%date: October 22, 2015
%%date : December 24, 2014

\vskip 0.7cm \noindent{\bf Abstract.} Recently, Z.-W. Sun introduced two kinds of polynomials related to the Delannoy numbers, and
proved some supercongruences on sums involving those polynomials. We deduce new summation formulas for squares of those polynomials and
use them to prove that certain rational sums involving even powers of those polynomials are integers whenever they are evaluated at integers. This confirms
two conjectures of Z.-W. Sun. We also conjecture that many of these results have neat $q$-analogues.


\vskip 2mm
\noindent{\it MR Subject Classifications}: 11A07, 11B65, 05A10

\section{Introduction}
It is well known that, for any $m,n\geqslant 0$, the number
\begin{align*}
\sum_{k=0}^{n}{n\choose k}{m\choose k}2^k=\sum_{k=0}^{n}{n\choose k}{n+m-k\choose n}, %\label{eq:delannoy}
\end{align*}
called a {\it Delannoy number}, counts lattice paths from $(0,0)$ to $(m,n)$ in which only east $(1,0)$, north $(0,1)$, and northeast $(1,1)$
steps are allowed.
Recently, Z.-W. Sun \cite{Sun} introduced the following polynomials
\begin{align*}
d_n(x)&=\sum_{k=0}^n{n\choose k}{x\choose k}2^k,\\
s_n(x)&=\sum_{k=0}^n{n\choose k}{x\choose k}{x+k\choose k},
\end{align*}
and established some interesting supercongruences involving $d_n(x)$ or $s_n(x)$, such as
\begin{align}
\sum_{k=0}^{p-1}(2k+1)d_k(x)^2
&\equiv\begin{cases}
-x\pmod{p^2},&\text{if $x\equiv 0\pmod{p}$}, \\
x+1\pmod{p^2},&\text{if $x\equiv -1\pmod{p}$}, \\
0\pmod{p^2},&\text{otherwise},
\end{cases} \label{eq:sun-one} \\[5pt]
\sum_{k=0}^{p-1}(2k+1)s_k(x)^2
&\equiv 0\pmod{p^2},  \label{eq:sun-two}
\end{align}
where $p$ is an odd prime and $x$ is a $p$-adic integer.

Recall that a polynomial $P(x)$ in $x$ with real coefficients is called {\it integer-valued}, if $P(x)\in\mathbb{Z}$ for all $x\in\mathbb{Z}$.
In this paper, we shall prove the following generalizations of \eqref{eq:sun-one} and \eqref{eq:sun-two}, which were originally conjectured
by Z.-W. Sun (see \cite[Conjectures 6.1 and 6.12]{Sun}).

\begin{thm}\label{thm:main}
Let $m$ and $n$ be positive integers. Then all of
\begin{align*}
&\frac{x(x+1)}{2n^2}\sum_{k=0}^{n-1}(2k+1)d_k(x)^2,\quad \frac{1}{n}\sum_{k=0}^{n-1}(2k+1)d_k(x)^{2m},\quad
\frac{1}{n}\sum_{k=0}^{n-1}(-1)^k(2k+1)d_k(x)^{2m}, \\
&\ \qquad{}\frac{1}{n^2}\sum_{k=0}^{n-1}(2k+1)s_k(x)^2,\quad \frac{1}{n}\sum_{k=0}^{n-1}(2k+1)s_k(x)^{2m},\quad
\frac{1}{n}\sum_{k=0}^{n-1}(-1)^k(2k+1)s_k(x)^{2m}
\end{align*}
are integer-valued.
\end{thm}

We shall also prove the following result, which will play an important role in our proof of Theorem \ref{thm:main}.
\begin{thm}\label{thm:second}
Let $m$ and $n$ be positive integers and let $j$, $k$ be non-negative integers. Then
\begin{align*}
\frac{(n-k)(k+1)}{n}{n+k\choose 2k}{m+1\choose k+1}{m+k\choose k+1}
\end{align*}
and
\begin{align*}
\frac{1}{k+1}{n-1\choose k}{n+k\choose k}{2k\choose j+k}{m+k\choose 2k}{m\choose j}{m+j\choose j}
\end{align*}
are integers.
\end{thm}

The paper is organized as follows. In the next section, we shall give a $q$-analogue of Theorem \ref{thm:second}.
In Section 3, we mainly give a single-sum expression for $d_n(x)^2$, a new expression for $s_n(x)^2$, and recall a recent divisibility result of Chen and Guo \cite{CG} concerning multi-variable Schmidt polynomials. The proof of Theorem \ref{thm:main} will be given in Section~4.
We propose some related open problems in the last section.

\section{A $q$-analogue of Theorem \ref{thm:second}}
Recall that the {\it $q$-binomial coefficients}
are defined by
\begin{align*}
{n\brack k}
=\begin{cases}
\displaystyle \prod_{i=1}^{k}\frac{1-q^{n-k+i}}{1-q^i},
&\text{if $0\leqslant k\leqslant n,$} \\[10pt]
0, &\text{otherwise.}
\end{cases}
\end{align*}
The following is our announced strengthening of Theorem \ref{thm:second}.

\begin{thm}\label{thm:q-analog}
Let $m$ and $n$ be positive integers and let $j$, $k$ be non-negative integers. Then
\begin{align}
\frac{(1-q^{n-k})(1-q^{k+1})}{(1-q)(1-q^n)}{n+k\brack 2k}{m+1\brack k+1}{m+k\brack k+1} \label{eq:q-analog}
\end{align}
and
\begin{align}
\frac{1-q}{1-q^{k+1}}{n-1\brack k}{n+k\brack k}{2k\brack j+k}{m+k\brack 2k}{m\brack j}{m+j\brack j}  \label{eq:q-analog-2}
\end{align}
are polynomials in $q$ with non-negative integer coefficients.
\end{thm}

\noindent{\it Proof of Theorem {\rm\ref{thm:q-analog}.}}
It suffices to show that \eqref{eq:q-analog} and \eqref{eq:q-analog-2} are polynomials in $q$ with integer coefficients, since the proof of the non-negativity is exactly
the same as those in \cite{Guo,GK}.  We shall accomplish this by decomposing $q$-binomial coefficients into cyclotomic polynomials.

It is well known that
$$
q^n-1=\prod _{d\mid n} ^{}\Phi_d(q),
$$
where $\Phi_d(q)$ denotes the $d$-th cyclotomic polynomial in $q$.
For any real number $x$, let $\lfloor x\rfloor$ denote the largest integer less than or equal to $x$.  Then
$$
\frac{(1-q^{n-k})(1-q^{k+1})}{(1-q)(1-q^n)}{n+k\brack 2k}{m+1\brack k+1}{m+k\brack k+1}
=\prod _{d=2} ^{\max\{m+k,n+k\}}\Phi_d(q)^{e_d},
$$
with
\begin{align*}
e_d&
=\chi(d\mid n-k)+\chi(d\mid k+1)-\chi(d\mid n)
+\fl{\frac {n+k} {d}}
+\fl{\frac {m+1} {d}}
+\fl{\frac {m+k} {d}} \notag\\
&\quad{}
-\fl{\frac {n-k} {d}}
-\fl{\frac {2k} {d}}
-\fl{\frac {m-k} {d}}
-\fl{\frac {m-1} {d}}
-2\fl{\frac {k+1} {d}}, %\label{eq:chi}
\end{align*}
where $\chi(\mathcal S)=1$ if $\mathcal S$ is
true and $\chi(\mathcal S)=0$ otherwise.
The number $e_d$ is obviously non-negative, unless
$d\mid n$, $d\nmid n-k$ and $d\nmid k+1$.

So, let us assume that $d\mid n$, $d\nmid n-k$ and $d\nmid k+1$. Since one of $k$ and $k+1$ is even, we must have $d\geqslant 3$.
Let $\{x\}=x-\lfloor x\rfloor$ denote the fraction part of $x$.
We consider three cases: If $0<\{\frac{k}{d}\}<\frac{1}{2}$, then
\begin{align*}
\left\lfloor\frac{n+k}{d}\right\rfloor-\left\lfloor\frac{n-k}{d}\right\rfloor-\left\lfloor\frac{2k}{d}\right\rfloor
=\left\lfloor\frac{d+k}{d}\right\rfloor-\left\lfloor\frac{d-k}{d}\right\rfloor-\left\lfloor\frac{2k}{d}\right\rfloor=1.
\end{align*}
Namely, $e_d$ is non-negative. If $\{\frac{k}{d}\}\geqslant\frac{1}{2}$ and $\{\frac{m}{d}\}\geqslant \{\frac{k}{d}\}$, then noticing that $d\nmid k+1$,
we have $\{\frac{m-1}{d}\}\geqslant\frac{1}{2}-\frac{1}{3}>0$, and so
\begin{align*}
\left\lfloor\frac{m+k}{d}\right\rfloor-\left\lfloor\frac{m-1}{d}\right\rfloor-\left\lfloor\frac{k+1}{d}\right\rfloor
=1.
\end{align*}
That is, $e_d$ is also non-negative. If $\{\frac{k}{d}\}\geqslant\frac{1}{2}$ and $\{\frac{m}{d}\}<\{\frac{k}{d}\}$,
then $\{\frac{k+1}{d}\}>\{\frac{k}{d}\}\geqslant\frac{1}{2}$, and so $\{\frac{m+1}{d}\}<\{\frac{k+1}{d}\}$, i.e.,
\begin{align*}
\left\lfloor\frac{m+1}{d}\right\rfloor-\left\lfloor\frac{m-k}{d}\right\rfloor-\left\lfloor\frac{k+1}{d}\right\rfloor
=1.
\end{align*}
which means that $e_d$ is still non-negative.
This completes the proof of polynomiality of \eqref{eq:q-analog}.

Similarly, we have
\begin{align*}
\frac{1-q}{1-q^{k+1}}{n-1\brack k}{n+k\brack k}{2k\brack j+k}{m+k\brack 2k}{m\brack j}{m+j\brack j}
=\prod _{d=2} ^{\max\{n+k,m+k,m+j\}}\Phi_d(q)^{e_d},
\end{align*}
with
\begin{align*}
e_d&
=-\chi(d\mid k+1)
+\fl{\frac {n-1} {d}}
+\fl{\frac {n+k} {d}}
+\fl{\frac {m+k} {d}}
+\fl{\frac {m+j} {d}}
-\fl{\frac {n} {d}}
 \notag\\
&\quad{}
-\fl{\frac {n-k-1} {d}}
-\fl{\frac {j+k} {d}}
-\fl{\frac {k-j} {d}}
-\fl{\frac {m-k} {d}}
-2\fl{\frac {k} {d}}
-\fl{\frac {m-j} {d}}
-2\fl{\frac {j} {d}}. %\label{eq:chi}
\end{align*}
The number $e_d$ is obviously non-negative, unless
$d\mid k+1$.

Let $d\geqslant 2$ be a positive integer and $d\mid k+1$. It is clear that $\{\frac{k}{d}\}=\frac{d-1}{d}\geqslant \frac{1}{2}$.
\begin{itemize}
\item If $d\nmid n$, then $\fl{\frac {n-1} {d}}-\fl{\frac {n} {d}}=0$, and
\begin{align}
\fl{\frac {m+k} {d}}-\fl{\frac {m-k} {d}}-2\fl{\frac {k} {d}}\geqslant \fl{\frac {2k} {d}}-2\fl{\frac {k} {d}}=1. \label{eq:ineq-one}
\end{align}
\item If $d\mid n$, then $\fl{\frac {n-1} {d}}-\fl{\frac {n} {d}}=-1$ and the inequality \eqref{eq:ineq-one} still holds.
We consider three subcases:
\begin{itemize}
\item[(i)] For $\{\frac{m-k}{d}\}\geqslant \frac{2}{d}$, there holds
\begin{align*}
\fl{\frac {m+k} {d}}-\fl{\frac {m-k} {d}}-2\fl{\frac{k}{d}}=2. %\label{eq:ineq-one}
\end{align*}
\item[(ii)] For $\{\frac{m-k}{d}\}=\frac{1}{d}$, we have $d\mid m$. If $d\nmid j$, then
\begin{align*}
\fl{\frac {m+j} {d}}-\fl{\frac {m-j} {d}}-2\fl{\frac{j}{d}}\geqslant \fl{\frac {m} {d}}-\fl{\frac {m-j}{d}}-\fl{\frac{j}{d}}=1,
\end{align*}
while if $d\mid j$, then
\begin{align}
\fl{\frac {n+k} {d}}
-\fl{\frac {n-k} {d}}
-\fl{\frac {j+k} {d}}
-\fl{\frac {k-j} {d}}
=\fl{\frac {2k} {d}}
-2\fl{\frac {k} {d}}
=1.  \label{eq:inq-three}
\end{align}
\item[(ii)] For $\{\frac{m-k}{d}\}=0$, we have $d\mid m+1$. If $d\nmid j$, then
\begin{align*}
\fl{\frac {m+j} {d}}-\fl{\frac {m-j} {d}}-2\fl{\frac{j}{d}}\geqslant \fl{\frac {m+j} {d}}-\fl{\frac {m}{d}}-\fl{\frac{j}{d}}=1,
\end{align*}
while if $d\mid j$, then the inequality \eqref{eq:inq-three} holds again.
\end{itemize}
\end{itemize}
Above all, we have proved that $e_d\geqslant 0$ in any case.  This completes the proof of polynomiality of \eqref{eq:q-analog-2}.
\qed

\noindent{\it Remark.} It was pointed out by the referee that a slightly shorter proof of $e_d\geqslant 0$
can be given by noticing that we may assume that $0\leqslant n,k,m<d$.

\medskip
It is easily seen that Theorem~\ref{thm:second} follows from Theorem~\ref{thm:q-analog} by letting $q\to 1$.

\section{Some auxiliary results}

Z.-W. Sun \cite[(1.7)]{Sun} noticed that $d_n(-x-1)=(-1)^n d_n(x)$, which is also demonstrated by a formula in \cite[p.~31]{Slater}.
This encourages us to find the following identity for $d_n(x)^2$.
\begin{lem}Let $n$ be a non-negative integer. Then
\begin{align}
d_n(x)^2=\sum_{k=0}^{n}{n+k\choose 2k}{x\choose k}{x+k\choose k}4^k. \label{eq:dnsquare}
\end{align}
\end{lem}

\pf Denote the right-hand side of \eqref{eq:dnsquare} by $S_n(x)$. Applying the Zeilberger algorithm (see \cite{Koepf,PWZ}), we have
\begin{align}
&(n+2)d_{n+2}(x)=(2x+1)d_{n+1}(x)+(n+1)d_n(x), \label{eq:zeil-1} \\[5pt]
&(n+1)^2S_n(x)-(n^2+4n+4x^2+4x+5)\left(S_{n+1}(x)+S_{n+2}(x)\right)+(n+3)^2S_{n+3}(x)=0.  \notag %\label{eq:zeil-2}
\end{align}
It follows from \eqref{eq:zeil-1} that
\begin{align}
&(n+2)^2d_{n+2}(x)^2 =(2x+1)^2d_{n+1}(x)^2+(n+1)^2d_n(x)^2+2(n+1)(2x+1)d_{n+1}(x)d_n(x), \label{eq:zeil-3}\\
&(n+2)d_{n+2}(x)d_{n+1}(x)=(2x+1)d_{n+1}(x)^2+(n+1)d_{n+1}(x)d_n(x).  \label{eq:zeil-4}
\end{align}
Substituting \eqref{eq:zeil-3} twice into \eqref{eq:zeil-4}, and making some simplification, we immediately get
\begin{align*}
(n+1)^2d_n(x)^2-(n^2+4n+4x^2+4x+5)\left(d_{n+1}(x)^2+d_{n+2}(x)^2\right)+(n+3)^2d_{n+3}(x)^2=0.
\end{align*}
Namely, the polynomials $d_n(x)^2$ and $S_n(x)$ satisfy the same recurrence. It is easy to see that
$d_n(x)^2=S_n(x)$ holds for $n=0,1,2$. This completes the proof of \eqref{eq:dnsquare}.
\qed

\noindent{\it Remark.} The hypergeometric form of \eqref{eq:dnsquare} is as follows:
\begin{align}
_{2}F_1
\left[\begin{array}{c}
-n,\, -x\\
1
\end{array};2
\right]^2
={}_{4}F_3
\left[\begin{array}{c}
-n,\, -x,\, n+1,\, x+1\\
1,\, 1,\, \frac{1}{2}
\end{array};1
\right].  \label{eq:3f2-one}
\end{align}
Wadim Zudilin (personal communication) pointed out that \eqref{eq:3f2-one} is a special case of the following identity \cite[p.~80, (2.5.32)]{Slater}:
\begin{align*}
_{2}F_1
\left[\begin{array}{c}
a,\, b\\
c
\end{array};z
\right]
{}_{2}F_1
\left[\begin{array}{c}
a,\, c-b\\
c
\end{array};z
\right]
=(1-z)^{-a}{}_{4}F_3
\left[\begin{array}{c}
a,\, b,\, c-a,\, c-b\\
c,\,\frac{c}{2},\, \frac{c+1}{2},
\end{array};\frac{-z^2}{4(1-z)}
\right],  %\label{eq:3f2-two}
\end{align*}
by noticing the identity \cite[p.~31, (1.7.1.3)]{Slater}:
\begin{align*}
(1-z)^{-a}{} _{2}F_1
\left[\begin{array}{c}
a,\, b\\
c
\end{array};-\frac{z}{1-z}
\right]
=
{}_{2}F_1
\left[\begin{array}{c}
a,\, c-b\\
c
\end{array};z
\right].
\end{align*}
Besides, the polynomial $d_n(x)$ is a particular case of classical Meixner orthogonal polynomials
(see http://homepage.tudelft.nl/11r49/askey/ch1/par9/par9.html).

We also need the following new expression for $s_n(x)^2$, which is crucial in dealing with the last three polynomials
in Theorem \ref{thm:main}.
\begin{lem}Let $n$ be a non-negative integer. Then
\begin{align}
s_n(x)^2=\sum_{k=0}^{n}{n+k\choose 2k}{x\choose k}{x+k\choose k}\sum_{j=0}^{k}{2k\choose j+k}{x\choose j}{x+j\choose j}. \label{eq:sn-square}
\end{align}
\end{lem}
\pf From the Pfaff-Saalsch\"utz identity \cite[(1.4)]{Andrews}, we deduce that (see the proof of \cite[Lemma 4.2]{GZ})
\begin{align*}
{x\choose j}{x+j\choose j}{x\choose k}{x+k\choose k}=\sum_{r=k}^{j+k}{j+k\choose j}{k\choose r-j}{r\choose k}{x\choose r}{x+r\choose r}.
\end{align*}
Therefore, comparing the coefficients of ${x\choose r}{x+r\choose r}$, we see that \eqref{eq:sn-square} is equivalent to
\begin{align}
&\hskip -2mm \sum_{j=0}^{n}\sum_{k=0}^{n}{n\choose j}{n\choose k}{j+k\choose j}{k\choose r-j}{r\choose k} \nonumber \\
&=\sum_{j=0}^{n}\sum_{k=0}^{n}{n+k\choose 2k}{2k\choose j+k}{j+k\choose j}{k\choose r-j}{r\choose k}. \label{eq:double-sum}
\end{align}

Denote the left-hand side of \eqref{eq:double-sum} by $A_n$ and the right-hand side of \eqref{eq:double-sum} by $B_n$.
Then the multi-Zeilberger algorithm gives the following recurrences of order $3$:
\begin{align}
&\hskip -2mm (n+3)^2(2n-r+5)(2n-r+6)A_{n+3}-(12n^4-4n^3r+3n^2r^2+110n^3-17n^2r \notag\\
&{}+14nr^2+394n^2-18nr+17r^2+650n+r+414)A_{n+2}
+(12n^4+4n^3r+3n^2r^2   \notag\\
&{}+82n^3+31n^2r+10nr^2+226n^2+74nr+9r^2+294n+57r+150)A_{n+1} \notag\\
&{}-(n+1)^2(2n+r+2)(2n+r+3)A_n=0, \label{eq:rec-one}\\[5pt]
&\hskip -2mm (n+3)(2n+3)(2n-r+5)(2n-r+6)B_{n+3}-(2n+5)(4n^3+4n^2r+nr^2+30n^2    \notag\\
&{}+17nr+r^2+72n+17r+54)B_{n+2}-(2n+3)(4n^3-4n^2r+nr^2+18n^2-15nr   \notag\\
&{}+3r^2+24n-13r+10)B_{n+1}+(n+1)(2n+5)(2n+r+2)(2n+r+3)B_n=0.  \label{eq:rec-two}
\end{align}
By induction on $n$, we may deduce from \eqref{eq:rec-one} and \eqref{eq:rec-two} that the numbers $A_n$ and $B_n$ also satisfy
the same recurrence  of order $2$:
\begin{align*}
&\hskip -2mm (n+2)(2n-r+3)(2n-r+4)A_{n+2}-(2n+3)(4n^2+r^2+12n+r+10)A_{n+1}  \\
&{}+(n+1)(2n+r+2)(2n+r+3)A_n=0.
\end{align*}
Moreover, it is clear that $A_0=B_0$ and $A_1=B_1$ for any $r$. This proves that $A_n=B_n$ holds for all $n$.
\qed

\noindent{\it Remark.} If we apply the multi-Zeilberger algorithm to the right-hand side of \eqref{eq:sn-square} directly,
then we will obtain a much more complicated recurrence  of order $7$. This is why we turn to consider
the equivalent form \eqref{eq:double-sum} of the identity \eqref{eq:sn-square}.
\medskip

The following result can be easily proved by induction on $n$.
\begin{lem}Let $n$ and $k$ be non-negative integers with $k\leqslant n$. Then
\begin{align}
\sum_{m=k}^{n-1}(2m+1){m+k\choose 2k}=\frac{n(n-k)}{k+1}{n+k\choose 2k}. \label{eq:simple}
\end{align}
\end{lem}

Let
$$
S_{n}(x_0,\ldots,x_{n})=\sum_{k=0}^n {n+k \choose 2k}{2k\choose k} x_k.
$$
be the multi-variable Schmidt polynomials. In order to prove Theorem \ref{thm:main}, we also need the following result, which
is a special case of \cite[Theorem 1.1]{CG}.

\begin{lem}\label{lem:schmidt}
Let $m$ and $n$ positive integers and $\varepsilon=\pm1$. Then all the coefficients in
\begin{align*}
\sum_{k=0}^{n-1}\varepsilon^k (2k+1)S_{k}(x_0,\ldots,x_{k})^m
\end{align*}
are multiples of $n$.
\end{lem}

\section{Proof of Theorem \ref{thm:main}}
Applying the identities \eqref{eq:dnsquare} and \eqref{eq:simple}, we have
\begin{align*}
\sum_{m=0}^{n-1}(2m+1)d_m(x)^2
&=\sum_{m=0}^{n-1}(2m+1)\sum_{k=0}^{m}{m+k\choose 2k}{x\choose k}{x+k\choose k}4^k \\
&=\sum_{k=0}^{n-1}\frac{n(n-k)}{k+1}{n+k\choose 2k}{x\choose k}{x+k\choose k}4^k.
\end{align*}
Therefore,
\begin{align}
\frac{x(x+1)}{2n^2}\sum_{m=0}^{n-1}(2m+1)d_m(x)^2
=\sum_{k=0}^{n-1}\frac{(n-k)(k+1)}{2n}{n+k\choose 2k}{x+1\choose k+1}{x+k\choose k+1}4^k. \label{eq:xx+1}
\end{align}
We now assume that $x$ is a positive integer in \eqref{eq:xx+1}. Then by Theorem \ref{thm:second} the $k$-th summand in the right-hand side of \eqref{eq:xx+1}
is an integer for $k\geqslant 1$, and is equal to ${x+1\choose 2}$ for $k=0$. This proves the first polynomial in Theorem \ref{thm:main} is integer-valued.

Similarly, applying \eqref{eq:sn-square} and \eqref{eq:simple}, we have
\begin{align}
\frac{1}{n^2}\sum_{m=0}^{n-1}(2m+1)s_m(x)^2
%&=\sum_{m=0}^{n-1}(2m+1)\sum_{k=0}^{m}{m+k\choose 2k}{x\choose k}{x+k\choose k}\sum_{j=0}^{k}{2k\choose j+k}{x\choose j}{x+j\choose j} \\
&=\sum_{k=0}^{n-1}\frac{1}{k+1}{n-1\choose k}{n+k\choose k}{x+k\choose 2k}\sum_{j=0}^{k}{2k\choose j+k}{x\choose j}{x+j\choose j}, \label{eq:double-sum-two}
\end{align}
which by Theorem \ref{thm:second} is clearly integer-valued.

For any non-negative integer $k$, let
\begin{align*}
x_k &:={x+k\choose 2k}4^k, \\
y_k &:={x+k\choose 2k}\sum_{j=0}^{k}{2k\choose j+k}{x\choose j}{x+j\choose j}.
\end{align*}
Then the identities \eqref{eq:dnsquare} and \eqref{eq:sn-square} may be respectively rewritten as
\begin{align*}
d_n(x)^2=\sum_{k=0}^{n}{n+k\choose 2k}{2k\choose k}x_k, \\
s_n(x)^2=\sum_{k=0}^{n}{n+k\choose 2k}{2k\choose k}y_k.
\end{align*}
It is clear that the numbers $x_0,\ldots,x_n$ and $y_0,\ldots,y_n$ are all integers when $x$ is an integer. By Lemma \ref{lem:schmidt}, we see that the other \
four polynomials in Theorem \ref{thm:main} are also integer-valued.

\section{Concluding remarks and open problems}
A special case of a well-known ${}_3F_2$ transformation formula in \cite[p.~142]{AAR} gives:
\begin{align*}
{}_3F_2\left[\begin{array}{c} -n,\,-x,\,x+1 \\ 1,\,1\end{array};1\right]
=\frac{(x+1)_n}{n!}
{}_3F_2\left[\begin{array}{c} -n,\,-x,\,-x \\ 1,\,-x-n\end{array};1\right],
\end{align*}
i.e.,
\begin{align*}
s_n(x)=\sum_{k=0}^{n}{n\choose k}{x\choose k}{x+n-k\choose n}.
\end{align*}
Let $p\geqslant 5$ be an odd prime. Z.-W. Sun \cite[Conjecture 6.11]{Sun} also conjectured that
\begin{align}
\sum_{k=0}^{p-1}(2k+1)s_k(x)^2
\equiv
\begin{cases}
\displaystyle\frac{3}{4}\left(\frac{-1}{p}\right)p^2 \pmod{p^4}, &\text{if $\displaystyle x=-\frac{1}{2}$,}\\[15pt]
\displaystyle\frac{7}{9}\left(\frac{-3}{p}\right)p^2 \pmod{p^4}, &\text{if $\displaystyle x=-\frac{1}{3}$,}\\[15pt]
\displaystyle\frac{13}{16}\left(\frac{-2}{p}\right)p^2 \pmod{p^4}, &\text{if $\displaystyle x=-\frac{1}{4}$,}\\[15pt]
\displaystyle\frac{31}{36}\left(\frac{-1}{p}\right)p^2 \pmod{p^4}, &\text{if $\displaystyle x=-\frac{1}{6}$,}
\end{cases} \label{eq:sun-final}
\end{align}
where $(\frac{\cdot}{p})$ denotes the Legendre symbol modulo $p$.

It is easy to see that, for $0\leqslant k\leqslant p-1$,
$$
{p-1\choose k}{p+k\choose k}
\equiv (-1)^k  \pmod{p^2}.
$$
Hence, letting $n=p$ in \eqref{eq:double-sum-two} and applying Theorem \ref{thm:second}, we immediately obtain
\begin{thm}Let $p$ be a prime and $x$ a $p$-adic integer. Then
\begin{align}
\sum_{k=0}^{p-1}(2k+1)s_k(x)^2
\equiv
p^2\sum_{k=0}^{p-1}\frac{(-1)^k}{k+1}{x+k\choose 2k}\sum_{j=0}^{k}{2k\choose j+k}{x\choose j}{x+j\choose j} \pmod{p^4}. \label{eq:cong-final}
\end{align}
\end{thm}

We believe that the congruence \eqref{eq:cong-final} can be utilized to prove Sun's conjectural congruence \eqref{eq:sun-final}.
Unfortunately, we are unable to accomplish this work. We hope that the interested reader can continue working on this problem.

An identity similar to \eqref{eq:dnsquare} is Clausen's identity \cite{Clausen}:
\begin{equation}
_2F_1\left[\begin{array}{c} a,\,b \\ a+b+\frac{1}{2}\end{array};x\right]^2
={}_3F_2\left[\begin{array}{c} 2a,\,2b,\,a+b \\ 2a+2b,\,a+b+\frac{1}{2}\end{array};x\right]
,\quad |x|<1.  \label{eq:clausen}
\end{equation}
It is well known that Clausen's identity \eqref{eq:clausen} has three different $q$-analogues (see \cite[(8.8.17) and (III.22)]{GR} and \cite{Jackson40,Jackson41}).
It is natural to ask
\begin{prob}{\rm Is there a $q$-analogue of the identity \eqref{eq:dnsquare}?}
\end{prob}

Dziemia\'nczuk \cite{Dzie} considered weighted Delannoy numbers. The natural $q$-Delannoy numbers (see \cite[p.~30]{Dzie}) are
$$D_q(m,n):=\sum_{k=0}^{n}q^{k\choose 2}{n\brack k}{n+m-k\brack n}.$$
It is easy to see that
$$D_{q^{-1}}(m,n)=q^{-mn}\sum_{k=0}^{n}q^{k+1\choose 2}{n\brack k}{n+m-k\brack n}.$$
It seems that a possible $q$-analogue of the left-hand side of \eqref{eq:dnsquare} should be\break $q^{mn}D_q(m,n)D_{q^{-1}}(m,n)$
rather than $D_q(m,n)^2$. However, it is quite difficult
to find the corresponding $q$-analogue of the right-hand side of \eqref{eq:dnsquare}.

The following conjecture is a $q$-analogue of \eqref{eq:sun-one} in the case where $x$ is a positive integer.
\begin{conj}\label{conj:3}
Let $p$ be an odd prime and $m$ a positive integer. Then
\begin{align*}
&\hskip -3mm \sum_{k=0}^{p-1}\frac{1-q^{2k+1}}{1-q}D_q(m,k) D_{q^{-1}}(m,k) q^{-k}\\
&\equiv\begin{cases}
\displaystyle\frac{1-q^{-2m}}{1-q^2}q\pmod{[p]^2},&\text{if $m\equiv 0\pmod{p}$}, \\[10pt]
\displaystyle\frac{1-q^{2m+2}}{1-q^2}q\pmod{[p]^2},&\text{if $m\equiv -1\pmod{p}$}, \\[10pt]
0\pmod{[p]^2},&\text{otherwise},
\end{cases}
\end{align*}
where $[p]=1+q+\cdots+q^{p-1}$.
\end{conj}

Furthermore, a fascinating $q$-analogue of the first three expressions in Theorem \ref{thm:main} seems to be true.
\begin{conj}\label{conj:4}
Let $m$, $n$, and $r$ be positive integers. Then all of
\begin{align*}
&\sum_{k=0}^{n-1}\frac{(1-q^{m})(1-q^{m+1})(1-q^{2k+1})}{(1-q^2)(1-q^n)^2}D_q(m,k)D_{q^{-1}}(m,k)q^{-k}, \\[5pt]
&\sum_{k=0}^{n-1}\frac{1-q^{2k+1}}{1-q^n}D_q(m,k)^r D_{q^{-1}}(m,k)^r q^{-k}, \\[5pt]
&\sum_{k=0}^{n-1}(-1)^{n-k-1}\frac{1-q^{2k+1}}{1-q^n}D_q(m,k)^r D_{q^{-1}}(m,k)^r q^{k\choose 2}
&\end{align*}
are Laurent polynomials in $q$ with non-negative integer coefficients.
\end{conj}

To prove Conjectures \ref{conj:3} and \ref{conj:4}, perhaps we need to give a single-sum expression for\break $D_q(m,n)D_{q^{-1}}(m,n)$ and
a $q$-analogue of Lemma \ref{lem:schmidt}. The latter is relatively easy, while the former is rather difficult because our proofs of \eqref{eq:dnsquare}
cannot be extended to the $q$-analogue case directly. By the way, we did not find any $q$-analogue of the other three polynomials in Theorem~\ref{thm:main}.

Finally, based on numerical calculations, we propose the following conjecture.
\begin{conj}Let $m$ and $n$ be positive integers. Then both
$$
\frac{1}{n}\sum_{k=0}^{n-1}(2k+1)d_k(x)^m s_k(x)^m\quad\text{and}\quad \frac{1}{n}\sum_{k=0}^{n-1}(-1)^k(2k+1)d_k(x)^m s_k(x)^m
$$
are integer-valued.
\end{conj}

\vskip 5mm \noindent{\bf Acknowledgments.} This work was partially supported by the National Natural Science Foundation of China (grant no. 11371144)
and the Qing Lan Project of Jiangsu Province.
The author would like to thank the editor and the anonymous referee for their helpful comments on a previous version of this paper.

\end{document}